\documentclass[11pt]{article}
\usepackage{graphics,amsmath,amssymb}
\usepackage{latexsym}
\usepackage{epsfig}
\usepackage{hyperref}
\usepackage{xcolor}
\usepackage{tikz}
\usetikzlibrary{calc,patterns,decorations.pathmorphing,decorations.markings,arrows}
\usepackage{array}   

\def\a{\alpha}
\def\be{\beta}

\def\epsilon{\varepsilon}
\def\ga{\gamma}

\def\la{\lambda}
\def\La{\Lambda}
\def\phi{\varphi}
\def\si{\sigma}

\def\om{\omega}

\newtheorem{theorem}{Theorem}[section]
\newtheorem{lemma}[theorem]{Lemma}

\newtheorem{corollary}[theorem]{Corollary}

\newtheorem{proposition}[theorem]{Proposition}
\newtheorem{remark}[theorem]{Remark}

\def\Z{{\mathbb Z}}
\def\N{{\mathbb N}}
\def\C{{\mathbb C}}
\def\R{{\mathbb R}}
\def\d{\,\hbox{\rm d}}

\newenvironment{Proof}{\removelastskip\par\medskip
\noindent{\em Proof.} \rm}{\penalty-20\null$\square$\par\medbreak}

\newenvironment{Proofy}{\removelastskip\par\medskip
\noindent{\em Proof of Theorem~\ref{th:inv.ingham}.} \rm}{\penalty-20\null$\square$\par\medbreak}

\newenvironment{Proofv}{\removelastskip\par\medskip
\noindent{\em Proof of Theorem~\ref{th:dir.ingham}.} \rm}{\penalty-20\null$\square$\par\medbreak}

\newenvironment{Proof1}{\removelastskip\par\medskip
\noindent{\em Proof of Theorem~\ref{th:diringham}.} \rm}{\penalty-20\null$\square$\par\medbreak}

\newenvironment{Proof2}{\removelastskip\par\medskip
\noindent{\em Proof of Proposition~\ref{pr:haraux-inv}.} \rm}{\penalty-20\null$\square$\par\medbreak}

\newenvironment{Proof3}{\removelastskip\par\medskip
\noindent{\em Proof of Theorem~\ref{th:inv.ingham1}.} \rm}{\penalty-20\null$\square$\par\medbreak}

\title{\bf Inverse Ingham type inequalities
\\ for the Burgers model}

\author{Paola Loreti
\thanks{Dipartimento di Scienze di Base e Applicate per l'Ingegneria,
Sapienza Universit\`a di Roma,
Via Antonio Scarpa 16, 00161 Roma (Italy); e-mail: $<$paola.loreti@uniroma1.it$>$ }
\and Daniela Sforza
\thanks{Dipartimento di Scienze di Base e Applicate per l'Ingegneria, 
Sapienza Universit\`a di Roma,
Via Antonio Scarpa 16, 00161 Roma (Italy); e-mail: $<$daniela.sforza@uniroma1.it$>$ }}

\begin{document}
%\date{}

\maketitle

\begin{abstract}
Viscoelastic materials have the properties both of elasticity and viscosity. In a previous work we investigate glass
relaxation in the framework of viscoelasticity. Here
we consider the Burgers  model, a first but
 meaningful step in the general analysis, showing a
reachability theorem thanks the analysis of the gap between eigenvalues and the representation of the solution
in Fourier series.
\end{abstract}



\section{Introduction}
The aim of this paper is to show 
exact reachability in time $T_0$ 
of Burgers equation. Burgers equation concerns  
 viscoelastic materials having the properties both of elasticity and viscosity. Also it can be seen
 as  first step to study relaxation of glasses 
 through approximation of the stretched exponential
 function with Prony series.
 As references on the topic see \cite{MM,MZ,MPVP}.
 The equation of viscoelasticity of Burgers models
is solved using  the Fourier analysis in \cite{LoretiSforza5}, indeed 
 in \cite{LoretiSforza5} we are able to give  a detailed
 spectral analysis and, by using Fourier expansions,
the solution of the Burgers equation.
Thanks to the previous results, in this paper we
state and solve for sufficiently large time the reachability
problem.  Theorem \ref{th:inv.ingham} , the inverse observability inequality
is shown  by using Ingham type approach, see \cite{Ing} and also \cite{KL2}.
This together with Theorem \ref{th:dir.ingham}, allow us to obtain the result
It remains an open question to  find the optimal time.
\begin{theorem}
Let $b_1,b_2,r_1,r_2>0$ with $\frac32(b_1+b_2)<r_1+r_2$ and $\frac{b_1}{r_1}+\frac{b_2}{r_2}=1$. 

Then, there exists $T_0>0$ such that for  $T>T_0$ and
$
(u_{0},u_{1})\in  L^{2}(0,\pi)\times H^{-1}(0,\pi)
$,
 there exists $f\in L^2(0,T)$ such that the weak solution
$u$ of the equation
\begin{equation*}
\displaystyle
u_{tt}(t,x) = u_{xx}(t,x)-b_1\int_0^t\ e^{-r_1(t-s)} u_{xx}(s,x)ds-b_2\int_0^t\ e^{-r_2(t-s)} u_{xx}(s,x)ds,
\quad t\in (0,T),\,\,\, x\in(0,\pi),
\end{equation*}
with boundary conditions
\begin{equation*}
u(t,0)=0,\quad u(t,\pi)=f(t),\qquad t\in (0,T),
\end{equation*}
and null initial values
\begin{equation*}
u(0,x)=u_{t}(0,x)=0,\qquad  x\in(0,\pi),
\end{equation*}
verifies the final conditions
\begin{equation*}
u(T,x)=u_{0}(x),\quad u_{t}(T,x)=u_{1}(x),
\quad x\in(0,\pi).
\end{equation*}
\end{theorem}
\section{Preliminaries}
Throughout the paper, we adopt the convention to write $F\asymp G$ if there exist two positive constants $c_1$ and $c_2$ such that
$c_1F\le G\le c_2F$.

Let $\Omega\subset\R^N$,  $N\ge1$, be a bounded open set with sufficiently smooth boundary $\partial\Omega$. 
We study the integro-differential equation %\eqref{eq:problem0}  can be 
written in the form
\begin{equation}\label{eq:problem01}
u_{tt} =\gamma^2\Delta u-b_1\int_0^t\ e^{-r_1(t-s)}\gamma^2\Delta u(s)ds-b_2\int_0^t\ e^{-r_2(t-s)}\gamma^2\Delta u(s)ds
\,,
\end{equation}
\begin{equation}
u(t,x)=0\qquad   t\ge0, \,\, x\in\partial\Omega,
\end{equation}
where $\gamma>0$ and $b_1, r_1,b_2, r_2$ are positive constants  satisfying
the conditions
\begin{equation}\label{cond:burgers}
r_1+ r_2-  b_1-  b_2>0,
\qquad
\frac{b_1}{r_1}+\frac{b_2}{r_2}=1\,.
\end{equation}
The assumptions \eqref{cond:burgers} on the integral kernels follow from the Burgers model.
%The constants $b_i$ and $r_i$ are defined in \eqref{eq:def-b} and \eqref{eq:Findley} respectively.

We can rewrite the integro-differential equation \eqref{eq:problem01}
in an abstract version. Let 
$
H= L^2(\Omega)
$
be endowed with the usual scalar product and norm 
\begin{equation*}
\langle u,v\rangle:=\int_\Omega u(x)v(x)\ dx\,,
\quad
\|u\|_H:=\left(\int_\Omega |u(x)|^{2}\ dx\right)^{1/2}
\quad
u,v\in H\,.
\end{equation*}
We define the operator $L:D(L)\subset H\to H$  by
\begin{equation}\label{eq:operatorL}
\begin{array}{l}
D(L)=H^2(\Omega)\cap H_0^1(\Omega) \\
\\
L u=\displaystyle -\gamma^2\Delta u\qquad u\in D(L)\,.
\end{array}
\end{equation}
It is well known that $L$ is a self-adjoint positive
 operator on $H$ with dense domain $D(L)$. The spectrum of $L$ is composed of an increasing sequence $\{\la_n\}_{n\ge1}$ of positive eigenvalues with
$\la_n\to\infty$ and there exists an orthonormal basis $\{e_n\}_{n\ge1}$ of $L^2(\Omega)$ consisting of the corresponding eigenvectors. Moreover, we assume that the eigenvalues $\la_n$ are all distinct numbers. 
 
%We denote by $\{\la_n\}_{n\ge1}$ a strictly increasing sequence  of  eigenvalues for the operator $L$ with
%$\la_n>0$ and $\la_n\to\infty$ and we assume that the sequence of the corresponding eigenvectors $\{e_n\}_{n\ge1}$ constitutes a Hilbert basis for $H$.

Recalling that $b_i, r_i>0$, $i=1,2$, %are defined in \eqref{eq:def-b} and \eqref{eq:Findley} and 
satisfy the condition $\frac{b_1}{r_1}+\frac{b_2}{r_2}=1$,
we consider the following Cauchy problem:
\begin{equation}\label{eq:system}
\begin{cases}
\displaystyle 
u''(t) +Lu(t)-b_1\int_0^t\ e^{-r_1(t-s)}L u(s)ds-b_2\int_0^t\ e^{-r_2(t-s)}L u(s)ds= 0
\hskip1cm  t\ge 0,
\\
u(0)=u_{0}\,,\quad u'(0)=u_{1}
\,.
\end{cases}
\end{equation}
For any couple $u_{0}\in D(\sqrt{L})$ and $u_{1}\in H$ of initial data we can give  an expansion by means of the eigenvectors $e_n$, that is
\begin{equation}\label{eq:v0}
\begin{split}
& u_{0}=\sum_{n=1}^{\infty}u_{0n}e_n\,,\qquad\quad u_{0n}=
\langle u_{0},e_n\rangle \,,
\quad
\|u_0\|^2_{D(\sqrt{L})}=\sum_{n=1}^{\infty}u_{0n}^2\lambda_n\,,
\\
& u_{1}=\sum_{n=1}^{\infty}u_{1n}e_n\,,\qquad\quad u_{1n}=\langle u_{1},e_n\rangle \,,\quad
\|u_1\|^2_{H}=\sum_{n=1}^{\infty}u_{1n}^2\,.
\end{split} 
\end{equation} 
\begin{theorem}\label{th:repres}
Given $u_{0}=\sum_{n=1}^{\infty}u_{0n}e_n\in D(\sqrt{L})$ and $u_{1}=\sum_{n=1}^{\infty}u_{1n}e_n\in H$,
for any $n\in\N$ we define the numbers 
\begin{equation}\label{eq:lambda2}
\omega_{n}=
\sqrt{\la_{n}}+i{b_1+b_2\over 2}
+O\Big({1\over{\sqrt{\la_{n}}}}\Big),
\end{equation}\begin{equation}\label{eq:lambda1}
\rho_{n}=b_1+b_2-r_1-r_2
+O\Big({1\over{\la_{n}}}\Big),
\end{equation}
\begin{equation}\label{eq:asy_C}
C_{n}={u_{0n}\over 2}
-\frac{i}{4}
\big((b_1+b_2) u_{0n}+2u_{1n}\big)\frac{1}{\sqrt{\lambda_n}}
+(u_{0n}+u_{1n})O\Big({1\over{\la_{n}}}\Big)
\,,
\end{equation}
\begin{equation}\label{eq:asy_R1}
R_{1,n}={r_1r_2u_{1n}\over(r_1+r_2-b_1-b_2)\la_{n}}
+( u_{0n}+u_{1n})O\Big({1\over{\la_n^{2}}}\Big),
\end{equation}
\begin{equation}\label{eq:asy_R2}
R_{2,n}={(b_1+b_2-r_1)(b_1+b_2-r_2)\big(u_{0n}(b_1+b_2-r_1-r_2)+u_{1n}\big)\over(b_1+b_2-r_1-r_2)\la_{n}}
+( u_{0n}+u_{1n})O\Big({1\over{\la_n^{2}}}\Big).
\end{equation}
Then, the series
\begin{equation}\label{eq:u1}
u(t) =\sum_{n=1}^{\infty}\Big(
C_{n}e^{i\omega_{n} t}+\overline{C_{n}}e^{-i\overline{\omega_{n}}t}+R_{1,n}+R_{2,n}e^{\rho_{n} t}\Big)e_n
\end{equation}
is the solution of problem \ref{eq:system}.

Moreover, for some constant $M>0$ we have
\begin{equation*}
\vert R_{1,n}\vert+\vert R_{2,n}\vert\le\frac{M}{\sqrt{\lambda_n}}|C_{n}| \qquad \forall n\in\N, 
\end{equation*}
\begin{equation}\label{eq:ine101bis}
\sum_{ n= 1}^{\infty}\la_{n}|C_{n}|^2
\asymp
\|u_0\|^2_{D(\sqrt{L})}+\|u_1\|^2_{H}
\,.
\end{equation}
\end{theorem}

\section{An Ingham type inequality}

Our goal is to prove an inverse inequality %and a direct inequality 
for the function 
\begin{equation}\label{eq:vsum1}
u(t) =\sum_{n=1}^{\infty}\Big(
C_{n}e^{i\omega_{n} t}+\overline{C_{n}}e^{-i\overline{\omega_{n}}t}+R_{1,n}+R_{2,n}e^{\rho_{n} t}\Big)\,,
\end{equation}
with $C_{n},\om_{n}\in\C$
and
$R_{1,n},R_{2,n},\rho_{n}\in\R$. 
Throughout this section we  assume that  there exist
$\gamma,\alpha_\omega>0$ and $\alpha_\rho<0$ such that 
\begin{equation}\label{eq:hom1}
\liminf_{n\to\infty}\big({\Re}\om_{n+1}-{\Re}\om_{n}\big)=\gamma\,,
\end{equation}
\begin{equation}\label{eq:hom2}
\lim_{n\to\infty}{\Im}\om_n=\alpha_\omega
\,,
\qquad
\lim_{n\to\infty}\rho_n=\alpha_\rho\,.
\end{equation}
Moreover, we suppose that  there exist
$\nu> 1/2$ and $M>0$
such that
\begin{equation}\label{eq:hom3}
|R_{1,n}|+|R_{2,n}|\le \frac{M}{n^{\nu}}|C_{n}|\qquad\forall\ n\in\N.
\end{equation}
\subsection{Some auxiliary results}
Here we list some results to be used later. 

For any $T>0$
we  consider the function defined by
\begin{equation}\label{eq:k}
g(t):=\left \{\begin{array}{l}
\displaystyle\sin \frac{\pi t}{T}\,\qquad\qquad \mbox{if}\,\, t\in\ [0,T]\,,\\
\\
0\,\qquad\qquad\quad\  \ \ \  \mbox{otherwise}\,.
\end{array}\right .
\end{equation}
Some  properties of $g$ are now in order.

\begin{lemma} \label{th:k}
Set
\begin{equation}\label{eqn:K}
G(w):=-\frac{T\pi}{T^2w^2-\pi^2}\,,\qquad w\in \C\,,
\end{equation}
the following properties hold.
\begin{itemize}
\item[(i)] 
 For any $w\in \C$ 
\begin{equation}\label{eqn:sinek1}
\int_{0}^{\infty} g(t)e^{iw t}dt
= (1+e^{iw T})G(w)
%=2e^{iu \frac{T}{2}}\cos\big(u \frac{T}{2}\big)G(u)
\,.
\end{equation}
 \begin{equation}\label{eqn:sinek2}
\overline{G(w)}=G(\overline{w})\,,
\qquad
\big|G(w)\big|=\big|G(\overline{w})\big|\,,
\end{equation}
%\item[(ii)] 
%For any $z_i,w_i\in \C$, $i=1,2$, one has
%\begin{multline}\label{eq:sinek2biss}
%\int_{0}^{\infty} k(t)\Re(z_1e^{iw_1 t})\Re(z_2e^{iw_2t})dt
%\\
%=\frac12 \Re\Big(z_1z_2(1+e^{i(w_1+w_2) T})K(w_1+w_2)
%+z_1\overline{z_2}(1+e^{i(w_1-\overline{w_2}) T})K(w_1-\overline{w_2})\Big)\,.
%\end{multline}
\item[(ii)] 
Let $\sigma>0$ and $n\in\N$. Then,
for $\displaystyle T>\frac{2\pi}\sigma$  and $w\in\C$, $|w|\ge\sigma n$, we have
\begin{equation}\label{eq:sinek3}
\big|G(w)\big|\le 
%\frac{T}{\pi}\Big(\frac{2\pi}{\gamma T}\Big)^2
\frac{4\pi}{T\sigma^2(4n^2-1)}
\,.
\end{equation}
\end{itemize}
\end{lemma}

\begin{lemma}\label{le:gap}
%Assume that
%\begin{equation*}
%\liminf_{n\to\infty}\big(\Re\omega_{n+1}-\Re\omega_{n}\big)=\gamma>0\,.
%\end{equation*}
For any $\varepsilon\in (0,1)$ there exists $n_0\in\N$ such that
\begin{equation}\label{eq:Re_gap1}
|\Re\omega_n-\Re\omega_m|\ge \gamma\sqrt{1-\varepsilon} |n-m|\,,\qquad\forall n\,,m\ge n_0\,,
\end{equation}
\begin{equation}\label{eq:Re_gap2}
\Re\omega_n\ge \gamma\sqrt{1-\varepsilon}\ n\,,\qquad\forall n\ge n_0\,.
\end{equation}
\end{lemma}

\begin{lemma}\label{le:stimaK}
For any  $\varepsilon\in (0,1)$, $T>\frac{2\pi}{\gamma\sqrt{1-\varepsilon}}$ and  $a>0$ there exists 
$n_0\in\N$ such that 
\begin{equation}\label{eq:minus0}
a|G( \omega_{n})|
\le\frac{\pi\varepsilon}{T\gamma^2(1-\varepsilon)}
\qquad\forall n\ge n_0,
\end{equation}
\begin{equation}\label{eq:minus1}
a\sum_{n=n_0}^\infty\ |G( \omega_{n})|
\le\frac{\pi\varepsilon}{T\gamma^2(1-\varepsilon)}
\,.
\end{equation}
\end{lemma}
\begin{Proof}
Thanks to \ref{eq:Re_gap2} and  \ref{eq:sinek3} we observe that for $\varepsilon\in (0,1)$ and $T>\frac{2\pi}{\gamma\sqrt{1-\varepsilon}}$, fixed $a>0$, there exists $n_0\in\N$ such that \begin{equation*}
a|G( \omega_{n})|
\le\frac{\pi}{T\gamma^2(1-\varepsilon)}
\frac{4a}{4n^2-1}
\le\frac{\pi}{T\gamma^2(1-\varepsilon)}\varepsilon 
\qquad\forall n\ge n_0,
\end{equation*}
\begin{equation*}
a\sum_{n=n_0}^\infty\ |G( \omega_{n})|
\le\frac{\pi }{T\gamma^2(1-\varepsilon)}
\sum_{n=n_0}^\infty\
\frac{4 a}{4n^2-1}
\le\frac{\pi}{T\gamma^2(1-\varepsilon)}\varepsilon 
\,,
\end{equation*}
and hence we get our statement.
\end{Proof}
We can state the following result, for the proof see e.g. \cite{LoretiSforza1}.
\begin{proposition}\label{pr:Fn}
Let $g$ the weight function defined by \eqref{eq:k}.
Suppose that
\begin{equation*}
\liminf_{n\to\infty}\big(\Re \omega_{n+1}-\Re \omega_{n}\big)=\gamma>0
\end{equation*}
and $\{C_n\}$ is a complex number sequence  with
$\sum_{n=1}^\infty\ |C_{n}|^2<+\infty$.

Then for any  $\varepsilon\in (0,1)$ and $T>\frac{2\pi}{\gamma\sqrt{1-\varepsilon}}$ 
 there exists 
$n_0=n_0(\varepsilon)\in\N$ independent of $T$ and $C_n$ such that we have
\begin{multline}\label{eq:Fn2'}
\int_{0}^{\infty} g(t)\Big( \sum_{n=n_0}^\infty  C_{n}e^{i \omega_{n}t}+\overline{C_{n}}e^{-i \overline{\omega_n}t}\Big)^2\ d t
\\
\ge
2\pi T\sum_{n=n_0}^\infty\bigg( \frac{1}{\pi^2+4T^2(\Im \omega_{n})^2}
-
\frac{4}{T^2\gamma^2(1-\varepsilon)}\bigg) (1+e^{-2\Im \omega_{n} T})|C_{n}|^2
\,.
\end{multline}
%\begin{multline}\label{eq:Fn2''}
%\int_{0}^{\infty} k(t)\Big| \sum_{n=n_0}^\infty  F_{n}e^{i \sigma_{n}t}+\overline{F_{n}}e^{-i \overline{\sigma_n}t}\Big|^2\ d t
%\\
%\le
%2\pi T\sum_{n=n_0}^\infty\bigg( \frac{1}{\pi^2+4T^2(\Im \sigma_{n})^2}
%+
%\frac{4}{T^2\gamma^2}(1+\varepsilon)\bigg) (1+e^{-2\Im \sigma_{n} T})|F_{n}|^2
%\,.
%\end{multline}
\end{proposition}
\begin{proposition}\label{pr:invine}
Assume that  
%there exist
%$\gamma,\alpha_\omega>0$ and $\alpha_\rho<0$ such that 
\begin{equation}\label{eq:hom11}
\liminf_{n\to\infty}\big({\Re}\om_{n+1}-{\Re}\om_{n}\big)=\gamma>0\,,
\end{equation}
\begin{equation}\label{eq:hom21}
\lim_{n\to\infty}\rho_n=\alpha_\rho<0\,,
\end{equation}
and  there exist
$\nu> 1/2$ and $M>0$
such that
\begin{equation}\label{eq:hom31}
|R_{1,n}|+|R_{2,n}|\le \frac{M}{n^{\nu}}|C_{n}|\qquad\forall\ n\in\N.
\end{equation}
Then for any  $\varepsilon\in (0,1)$ and $T>\frac{2\pi}{\gamma\sqrt{1-\varepsilon}}$ 
 there exists 
$n_0=n_0(\varepsilon)\in\N$ independent of $T$ and $C_n$ such that we have
\begin{multline}\label{eq:ine10}
\int_{0}^{\infty} \Big(\sum_{n=n_0}^{\infty}
C_{n}e^{i\omega_{n} t}+\overline{C_{n}}e^{-i\overline{\omega_{n}}t}+R_{1,n}+R_{2,n}e^{\rho_{n} t}\Big)^2\ dt
\\
\ge
2\pi T\sum_{n=n_0}^\infty
\bigg( \frac{1}{\pi^2+4T^2(\Im \omega_{n})^2}-\frac{4(1+\varepsilon)}{T^2\gamma^2(1-\varepsilon)}\bigg) (1+e^{-2\Im \omega_{n} T})|C_{n}|^2
\,.
\end{multline}
\end{proposition}
\begin{Proof} We consider the weight function $g$ defined by \eqref{eq:k}.
\begin{multline*}
\int_{0}^{\infty} g(t)\Big(\sum_{n=n_0}^{\infty}
C_{n}e^{i\omega_{n} t}+\overline{C_{n}}e^{-i\overline{\omega_{n}}t}+R_{1,n}+R_{2,n}e^{\rho_{n} t}\Big)^2\ dt
\\ 
=
\int_{0}^{\infty} g(t)\Big( \sum_{n=n_0}^\infty  C_{n}e^{i \omega_{n}t}+\overline{C_{n}}e^{-i \overline{\omega_n}t}\Big)^2\ d t
+2\sum_{n,m=n_0}^{\infty}\int_{0}^{\infty} g(t)
\Big(C_{n}e^{i\omega_{n} t}+\overline{C_{n}}e^{-i\overline{\omega_{n}}t}\Big)\Big(R_{1,m}+R_{2,m}e^{\rho_{m} t}\Big)\ dt
\\
+\int_{0}^{\infty} g(t)\Big(\sum_{n=n_0}^{\infty}
R_{1,n}+R_{2,n}e^{\rho_{n} t}\Big)^2\ dt\,.
\end{multline*}
Since $g(t)$ is positive we have
\begin{multline*}
\int_{0}^{\infty} g(t)\Big(\sum_{n=n_0}^{\infty}
C_{n}e^{i\omega_{n} t}+\overline{C_{n}}e^{-i\overline{\omega_{n}}t}+R_{1,n}+R_{2,n}e^{\rho_{n} t}\Big)^2\ dt
\ge
\int_{0}^{\infty} g(t)\Big( \sum_{n=n_0}^\infty  C_{n}e^{i \omega_{n}t}+\overline{C_{n}}e^{-i \overline{\omega_n}t}\Big)^2\ d t
\\
+2\sum_{n,m=n_0}^{\infty}\int_{0}^{\infty} g(t)
\Big(C_{n}e^{i\omega_{n} t}+\overline{C_{n}}e^{-i\overline{\omega_{n}}t}\Big)\Big(R_{1,m}+R_{2,m}e^{\rho_{m} t}\Big)\ dt
\,.
\end{multline*}
We estimate the first term on the right-hand side of the previous inequality by means of Proposition \ref{pr:Fn}, and hence, thanks to \eqref{eq:Fn2'} we deduce
\begin{multline}\label{eqn:sum-s}
\int_{0}^{\infty} g(t)\Big(\sum_{n=n_0}^{\infty}
C_{n}e^{i\omega_{n} t}+\overline{C_{n}}e^{-i\overline{\omega_{n}}t}+R_{1,n}+R_{2,n}e^{\rho_{n} t}\Big)^2\ dt
\\
\ge
2\pi T\sum_{n=n_0}^\infty\bigg( \frac{1}{\pi^2+4T^2(\Im \omega_{n})^2}
-
\frac{4}{T^2\gamma^2(1-\varepsilon)}\bigg) (1+e^{-2\Im \omega_{n} T})|C_{n}|^2
\\
+2\sum_{n,m=n_0}^{\infty}\int_{0}^{\infty} g(t)
\Big(C_{n}e^{i\omega_{n} t}+\overline{C_{n}}e^{-i\overline{\omega_{n}}t}\Big)\Big(R_{1,m}+R_{2,m}e^{\rho_{m} t}\Big)\ dt
\,.
\end{multline}
 As regards the second term, we note that
\begin{multline}\label{eq:2term}
\sum_{n,m=n_0}^{\infty}\int_{0}^{\infty} g(t)
\Big(C_{n}e^{i\omega_{n} t}+\overline{C_{n}}e^{-i\overline{\omega_{n}}t}\Big)\Big(R_{1,m}+R_{2,m}e^{\rho_{m} t}\Big)\ dt
\\
=
\sum_{n,m=n_0}^{\infty}R_{1,m}\Re \big[C_n(1+e^{ i\om_n T})G( \om_n)\big]
+\sum_{n,m=n_0}^{\infty}R_{2,m}\Re \big[C_n(1+e^{( i\om_n+\rho_{m}) T})G( \om_n-i\rho_{m})\big]
\\
\ge
-\sum_{n,m=n_0}^{\infty}|R_{1,m}| |C_n|(1+e^{ -\Im\om_n T})\big|G( \om_n)\big|
-\sum_{n,m=n_0}^{\infty}|R_{2,m}| |C_n|(1+e^{(\rho_{m}-\Im\om_n) T})\big|G( \om_n-i\rho_{m})\big|
\,.
\end{multline}
Thanks to \eqref{eq:hom3} we have
\begin{multline*}
\sum_{n,m=n_0}^{\infty}|R_{1,m}||C_n|(1+e^{-\Im\om_n T}) \big|G( \om_n)\big| 
\le  M
\sum_{n,m=n_0}^{\infty}|C_m|\frac{|C_n|}{m^\nu} (1+e^{-\Im\om_nT}) \big|G( \om_n)\big| 
\\
\le
M\sum_{m=n_0}^{\infty}|C_m|^2\sum_{n=n_0}^{\infty}\big|G( \om_n)\big| 
+M\sum_{n=n_0}^{\infty}|C_n|^2(1+e^{-2\Im\om_nT})\big|G( \om_n)\big|\sum_{m=n_0}^{\infty}\frac{1 }{m^{2\nu}}
\,.
\end{multline*}
By using respectively \eqref{eq:minus1} with $a=M$ and \eqref{eq:minus0} with $a=M\sum_{m=1}^{\infty}\frac{1 }{m^{2\nu}}$ we obtain
\begin{equation*}
M\sum_{m=n_0}^{\infty}|C_m|^2\sum_{n=n_0}^{\infty}\big|G( \om_n)\big| 
\le
\frac{\pi\varepsilon}{T\gamma^2(1-\varepsilon)}\sum_{n=n_0}^{\infty}|C_n|^2(1+e^{-2\Im\om_nT}),
\end{equation*}
\begin{equation*}
M\sum_{n=n_0}^{\infty}|C_n|^2(1+e^{-2\Im\om_nT})\big|G( \om_n)\big|\sum_{m=n_0}^{\infty}\frac{1 }{m^{2\nu}}
\le
\frac{\pi\varepsilon}{T\gamma^2(1-\varepsilon)}\sum_{n=n_0}^{\infty}|C_n|^2(1+e^{-2\Im\om_nT}),
\end{equation*}
whence
\begin{equation}\label{eqn:mistisin0}
\sum_{n,m=n_0}^{\infty}|R_{1,m}||C_n|(1+e^{-\Im\om_n T}) \big|G( \om_n)\big| 
\le
\frac{2\pi\varepsilon}{T\gamma^2(1-\varepsilon)}\sum_{n=n_0}^{\infty}|C_n|^2(1+e^{-2\Im\om_nT})
\,.
\end{equation}
In a similar way, taking also into account  \eqref{eq:hom2} and $\alpha_\rho<0$, one can get
\begin{equation}\label{eqn:mistisin}
\sum_{n,m=n_0}^{\infty}|R_{2,m}||C_n|(1+e^{(\rho_{m}-\Im\om_n) T}) \big|G( \om_n-i\rho_{m})\big| 
\le
\frac{2\pi\varepsilon}{T\gamma^2(1-\varepsilon)}\sum_{n=n_0}^{\infty}|C_n|^2(1+e^{-2\Im\om_nT})
\,.
\end{equation}
Plugging \eqref{eqn:mistisin0} and \eqref{eqn:mistisin} into \eqref{eq:2term} we obtain
\begin{equation*}
\sum_{n,m=n_0}^{\infty}\int_{0}^{\infty} g(t)
\Big(C_{n}e^{i\omega_{n} t}+\overline{C_{n}}e^{-i\overline{\omega_{n}}t}\Big)\Big(R_{1,m}+R_{2,m}e^{\rho_{m} t}\Big)\ dt
\ge
-\frac{4\pi\varepsilon}{T\gamma^2(1-\varepsilon)}\sum_{n=n_0}^{\infty}|C_n|^2(1+e^{-2\Im\om_nT})
\,.
\end{equation*}
By \eqref{eqn:sum-s} we get
\begin{multline*}
\int_{0}^{\infty} g(t)\Big(\sum_{n=n_0}^{\infty}
C_{n}e^{i\omega_{n} t}+\overline{C_{n}}e^{-i\overline{\omega_{n}}t}+R_{1,n}+R_{2,n}e^{\rho_{n} t}\Big)^2\ dt
\\
\ge
2\pi T\sum_{n=n_0}^\infty
\bigg( \frac{1}{\pi^2+4T^2(\Im \omega_{n})^2}-\frac{4(1+\varepsilon)}{T^2\gamma^2(1-\varepsilon)}\bigg) (1+e^{-2\Im \omega_{n} T})|C_{n}|^2
\,.
\end{multline*}
In conclusion, taking into account the definition of $g$ \eqref{eq:ine10} follows from the above inequality.
\end{Proof}
To make meaningful the inequality \eqref{eq:ine10} we need to discuss the right-hand side.
\begin{theorem}
Under the assumptions of Proposition \ref{pr:invine}, if
$\gamma>4\alpha_\omega$ and 
$T>\frac{2\pi}{\sqrt{\gamma^2(1-\varepsilon)-16\alpha_\omega^2(1+\varepsilon)}}$, 
then there exists 
$n_0=n_0(\varepsilon)\in\N$ independent of $T$ and $C_n$ such that we have
\begin{multline}\label{eq:ine101}
\int_{0}^{T} \Big(\sum_{n=n_0}^{\infty}
C_{n}e^{i\omega_{n} t}+\overline{C_{n}}e^{-i\overline{\omega_{n}}t}+R_{1,n}+R_{2,n}e^{\rho_{n} t}\Big)^2\ dt
\\
\ge
2\pi T\bigg( \frac{1}{\pi^2+4T^2\alpha_\omega^2(1+\varepsilon)}
-\frac{4}{T^2\gamma^2(1-\varepsilon)}\bigg) \sum_{n=n_0}^\infty
(1+e^{-2\Im \omega_{n} T})|C_{n}|^2
\,.
\end{multline}
\end{theorem}
\begin{Proof}
We start by 
using \eqref{eq:hom2}
and replacing in \eqref{eq:ine10} $\varepsilon$ with $\varepsilon'\in (0,\frac{\varepsilon}{2-\varepsilon})$;  taking into account that $\frac{1+\varepsilon'}{1-\varepsilon'}<\frac1{1-\varepsilon}$, we have for $n\ge n_0$
\begin{equation*}
\frac{1}{\pi^2+4T^2(\Im \omega_{n})^2}-\frac{4(1+\varepsilon')}{T^2\gamma^2(1-\varepsilon')}
\ge
\frac{1}{\pi^2+4T^2\alpha_\omega^2(1+\varepsilon)}-\frac{4}{T^2\gamma^2(1-\varepsilon)}
\end{equation*}
The constant 
\begin{equation*}
\frac{1}{\pi^2+4T^2\alpha_\omega^2(1+\varepsilon)}
-\frac{4}{T^2\gamma^2(1-\varepsilon)}
\end{equation*}
is positive if
\begin{equation}\label{eq:const}
T^2\big[\gamma^2(1-\varepsilon)-16\alpha_\omega^2(1+\varepsilon)\big]
>4\pi^2\,.
\end{equation}
Since $\gamma>4\alpha_\omega$ we have $\gamma^2(1-\varepsilon)-16\alpha_\omega^2(1+\varepsilon)>0$
if $\varepsilon<\frac{\gamma^2-16\alpha_\omega^2}{\gamma^2+16\alpha_\omega^2}$.
If we assume the more restrictive condition $T>\frac{2\pi}{\sqrt{\gamma^2(1-\varepsilon)-16\alpha_\omega^2(1+\varepsilon)}}$ with respect to that
$T>\frac{2\pi}{\gamma\sqrt{1-\varepsilon}}$, then \ref{eq:const} holds true.
\end{Proof}
If we take the sequences
$\{\om_n\}_{n\in\Z}$ and
$\{\rho_n\}_{n\in\Z}$ composed by pairwise  distinct non null numbers
such that
$\rho_n\not=i\om_m$ for any $n\,,m\in\Z$ and
argue in a similar way as in \cite{LoretiSforza1,LoretiSforza3}, where the Haraux method \cite{Ha} is successfully applied, we get the following inverse inequality.
\begin{theorem}\label{th:inv.ingham}
For any 
$T>\frac{2\pi}{\sqrt{\gamma^2-16\alpha_\omega^2}}$ there exists a positive constant $C=C(T)$ such that
\begin{equation}\label{eq:ine101}
\int_{0}^{T} \Big(\sum_{n=1}^{\infty}
C_{n}e^{i\omega_{n} t}+\overline{C_{n}}e^{-i\overline{\omega_{n}}t}+R_{1,n}+R_{2,n}e^{\rho_{n} t}\Big)^2\ dt
\ge
C \sum_{n=1}^\infty
(1+e^{-2\Im \omega_{n} T})|C_{n}|^2
\,.
\end{equation}
\end{theorem}
%\textcolor{red}{rendiconti Parma per disuguaglianza diretta}
If we assume that $k: [0,\infty)\to [0,\infty)$ is a locally absolutely continuous
function such that $k(0)> 0$, $k'(t)\le 0$ for a.e. $t\ge 0$ and
$\int_0^{t} k (s)\ ds<1$, $t\ge 0$, then for the Cauchy problem 
\begin{equation}\label{eq:cauchyI}
\begin{cases}
\displaystyle
u_{tt}(t,x) - \triangle u(t,x)+\int_0^t\ k(t-s)
\triangle u(s,x)ds= 0\,,
\quad t\ge0,\,\, x\in \Omega,
\\
u(t,x)=0\qquad   t\ge0, \,\, x\in\partial\Omega,
\\
u(0,x)=u_{0}(x),\quad
u_t(0,x)=u_{1}(x),\qquad  x\in \Omega.
\end{cases}
\end{equation}
we can
adapt the argumentations done in \cite[Theorem 1.3]{LoretiSforza4} to obtain the following result.

\begin{theorem}\label{th:dir.ingham}
For $T>0$, there exists a constant $C_0>0$ depending on $T$ and $k$ such that for every $u_0\in
H^1_0(\Omega)$ and $u_1\in L^2(\Omega)$, denoted by  $u$ the weak solution  of \eqref{eq:cauchyI},
one can define $\partial_\nu u$ such that the following inequality holds
\begin{equation}\label{eq:dir.ingham}
\int_0^T\int_{\partial\Omega} |\partial_\nu u|^2d\sigma dt
\le C_0(\|u_0\|^2_{H^1_0(\Omega)}+\|u_1\|^2_{L^2(\Omega)})
\,.
\end{equation}
\end{theorem}
\section{Controllability}
\begin{theorem}
Let $b_1,b_2,r_1,r_2>0$ with $\frac32(b_1+b_2)<r_1+r_2$ and $\frac{b_1}{r_1}+\frac{b_2}{r_2}=1$. 

Then, there exists $T_0>0$ such that for  $T>T_0$  and
$
(u_{0},u_{1})\in  L^{2}(0,\pi)\times H^{-1}(0,\pi)
$,
 there exists $f\in L^2(0,T)$ such that the weak solution
$u$ of the equation
\begin{equation*}
\displaystyle
u_{tt}(t,x) = \gamma^2u_{xx}(t,x)-b_1\int_0^t\ e^{-r_1(t-s)} \gamma^2u_{xx}(s,x)ds-b_2\int_0^t\ e^{-r_2(t-s)}\gamma^2 u_{xx}(s,x)ds,
\quad t\in (0,T),\,\, x\in(0,\pi),
\end{equation*}
with boundary conditions
\begin{equation*}
u(t,0)=0,\quad u(t,\pi)=f(t),\qquad t\in (0,T),
\end{equation*}
and null initial values
\begin{equation*}
u(0,x)=u_{t}(0,x)=0,\qquad  x\in(0,\pi),
\end{equation*}
verifies the final conditions
\begin{equation*}
u(T,x)=u_{0}(x),\quad u_{t}(T,x)=u_{1}(x),
\quad x\in(0,\pi).
\end{equation*}
\end{theorem}
\begin{Proof}
%We can rewrite the equation in (\ref{eq:problem-usix}) as an abstract one of the type (\ref{eq:cauchy}).
To prove the statement we apply the Hilbert Uniqueness Method, following a similar strategy to that used in \cite{LoretiSforza1}.

Let $\Omega=(0,\pi )$ and 
$
X= L^2(0,\pi )
$
be endowed with the usual scalar product and norm 
$$
\|u\|:=\left(\int_0^\pi |u(x)|^{2}\ dx\right)^{1/2}\qquad
u\in L^2(0,\pi)\,.
$$
We consider the operator $L:D(L)\subset L^2(0,\pi )\to L^2(0,\pi )$ defined by 
\begin{equation}\label{eq:operatorL}
Lu=\displaystyle -\gamma^2u_{xx}
\qquad u\in D(L):=H^2(0,\pi )\cap H_0^1(0,\pi )\,.
\end{equation}
It is well known that $L$ is a self-adjoint positive
 operator on $L^2(0,\pi )$ with dense domain $D(L)$ and 
 $$D(\sqrt L)=H_0^1(0,\pi ).$$
 Moreover, the eigenvalues of  $L$ are $\gamma^2n^2$, $n\in\N$, and 
  the corresponding eigenvectors are given by $\sqrt{\frac2\pi}\sin( nx)$ that
constitute a Hilbert basis for $L^2(0,\pi )$.

We can apply our spectral analysis to the adjoint problem. 
First, we consider the adjoint equation given by
\begin{multline}\label{eq:adjoint1}
\displaystyle
z_{tt}(t,x) = \gamma^2z_{xx}(t,x)-b_1\int_t^T\ e^{-r_1(s-t)} \gamma^2z_{xx}(s,x)ds-b_2\int_t^T\ e^{-r_2(s-t)}\gamma^2 z_{xx}(s,x)ds,
\\
\ \ t\in (0,T),\,\, x\in(0,\pi),
\end{multline}
with the Dirichlet boundary condition
\begin{equation}\label{eq:bcond1}
z(t,0)=z(t,\pi)=0\qquad t\in (0,T),
\end{equation}
and  final data 
\begin{equation} \label{eq:final1}
z(T,\cdot)=z_0\in H^1_0(\Omega)\,,\qquad z_t(T,\cdot)=z_1\in L^2(\Omega)\,.
\end{equation}
%where $z_0\in H^1_0(\Omega)$ and $z_1\in L^2(\Omega)$.
The backward problem (\ref{eq:adjoint1})--(\ref{eq:final1}) is
equivalent to a Cauchy problem  with $u(t,x)=z(T-t,x)$.
Therefore we can apply the conclusions of the previous sections. First, we write the solution $z(t,x)$ of the adjoint problem  as a Fourier series.
Indeed, the solution $z$ of the adjoint problem can be written in the form 
\begin{equation*}
z(t,x)=\sum_{n=1}^{\infty}\Big(
C_{n}e^{i\omega_{n} (T-t)}+\overline{C_{n}}e^{-i\overline{\omega_{n}}(T-t)}+R_{1,n}+R_{2,n}e^{\rho_{n} (T-t)}\Big)\sin(nx)\qquad
(t,x)\in [0,T]\times [0,\pi]\,,
\end{equation*}
whence
\begin{equation*}
z_x(t,\pi)=\sum_{n=1}^{\infty}(-1)^n n\Big(
C_{n}e^{i\omega_{n} (T-t)}+\overline{C_{n}}e^{-i\overline{\omega_{n}}(T-t)}+R_{1,n}+R_{2,n}e^{\rho_{n} (T-t)}\Big)\qquad
t\in [0,T]\,.
\end{equation*}
We can apply theorems \ref{th:inv.ingham} and \ref{th:dir.ingham} to function $z_x(t,\pi)$.
%Therefore, thanks to inequalities  \eqref{eq:diringham} and \eqref{eq:inv.ingham}  theorem \ref{th:uniqueness} holds true. In addition,  
By estimates \eqref{eq:ine101},  \eqref{eq:ine101bis} and \eqref{eq:dir.ingham} we have that 
$$
\int_0^T |z_x(t,\pi)|^2\ dt
\asymp
\|z_0\|^2_{H^1_0(\Omega)}+\|z_1\|^2_{L^2(\Omega)}\,.
$$
%whence the space $F$ introduced at the end of section \ref{se:HUM} is $H^1_0(0,\pi)\times L^2(0,\pi)$. 
The proof is complete.
\end{Proof}

%\section{Conclusions}
%In this paper we have investigated glass relaxation models, starting by a well-known model in literature, see e.g. \cite{MM} and references therein. Due to the complexity of the problem we have approximated the stretched exponential relaxation by means of a Prony series. For a general Prony series  we have established some partial results concerning the spectral analysis of the problem. In particular, by induction on the number of the terms of the Prony series the integro-differential equation showing the viscoelastic properties of the glass relaxation has always a null eigenvalue and the sum of all its eigenvalues is given by minus the sum of 
%the exponents of the Prony series.  
%
%
%In order to give more accurate results, we simplified the problem by taking under consideration the Burgers model, where the Prony series consists of two decreasing exponential functions.
%In that case we  have been able to
%give a complete description of the oscillations of the material in its relaxation stage, when it shows viscoelastic features.
%In particular, our analysis has revealed that the accumulation point of the branch of the real eigenvalues $\rho_n$, see \eqref{eq:lambda1},  depends only on the Kelvin-Voigt unit $(E_2,\eta_2 )$, see Figure \ref{fig:Burgers}. Indeed, taking into account \eqref{eq:accp} and \eqref{eq:Findley},
%we have obtained
%\begin{equation*}
%b_1+b_2-r_1-r_2=-\frac{q_1}{q_2}=-\frac{E_2}{\eta_2}
%\,.
%\end{equation*}


\begin{thebibliography}{99}

\itemsep=\smallskipamount

%\bibitem{Al} F. Alabau-Boussouira
%{\em  A Two-Level Energy Method for Indirect Boundary Observability and Controllability of Weakly Coupled Hyperbolic Systems} SIAM J. Control Optim., {\bf 42} (2003), 871--906.


%\bibitem{Boltz} L. Boltzmann, {\em Zur Theorie der elastichen Nachwirkung}, Wiener Berichte, {\bf 70} (1874), 275-306.
%
%\bibitem{CCM} M. Cardona, R. V. Chamberlin, W. Marx, The history of the stretched exponential function, Ann.Phys., 16  (2007), 842.

%\bibitem{Ce} E. Ces\`aro, {\em Sur la convergence des s\'eries}, Nouvelles annales de math\'ematiques, {\bf 7} (1888), 49--59.
%(teorema a p.54)

%\bibitem{CN} B. D. Coleman, W. Noll, {\em Foundations of linear viscoelasticity}, Rev. Modern Phys., {\bf 33} (1961) 239--249. 
%
%\bibitem{D1} C. M. Dafermos, {\em Asymptotic stability in
%viscoelasticity}, Arch. Rational Mech. Anal., {\bf 37} (1970),
%297--308.
%
%\bibitem{D2} C. M. Dafermos,  {\em An abstract Volterra equation with
%applications to linear viscoelasticity}, J.
%Differential Equations, {\bf 7} (1970), 554--569.

%\bibitem{DE} M. Doi, S. F. Edwards,  Dynamics of concentrated polymer systems, Parts 1, 2 and 3, J. Chem. Soc. Faraday II  74 (1978), 1789--1832;
%Parts 4, J. Chem. Soc. Faraday II  75 (1979), 38--54.
%
%\bibitem{FLO} W. N. Findley, J. S. Lai, K. Onaran, {\em Creep and relaxation of nonlinear viscoplastic materials},
%North-Holland, New York (1976).

%\bibitem{GLS} G. Gripenberg, S. O. Londen, O. J. Staffans, {\em Volterra
%Integral and Functional Equations}, Encyclopedia Math. Applications, {\bf 34}
%(1990), Cambridge Univ. Press, Cambridge.

\bibitem{Ha} A.  Haraux, {\em S\'eries lacunaires et contr\^ole semi-interne des vibrations d'une plaque rectangulaire} J. Math. Pures Appl., {\bf 68} (1989), 
457--465.

\bibitem{Ing} A. E. Ingham, {\em Some trigonometrical inequalities with
applications to the theory of series}, Math. Z., {\bf 41} (1936), 367-379.


%\bibitem{K2} J. U. Kim,  {\em Control of a second-order integro-differential equation}, SIAM J. Control Optim., {\bf 31} (1993),  101--110.

%\bibitem{KR} V. Komornik, B. Rao,  {\em Boundary stabilization of compactly coupled wave equations} Asymptot. Anal. {\bf 14} (1997),  339--359.

%\bibitem{KL1}  V. Komornik, P. Loreti, {\em Ingham type theorems for
%vector-valued functions and observability of coupled linear system},
% SIAM J. Control Optim., {\bf 37} (1998), 461-485.

\bibitem{KL2}  V. Komornik, P. Loreti, {\em Fourier series in control theory},
Springer Monographs in Ma\-the\-ma\-tics (2005), Springer-Verlag, New York.

%\bibitem{Las} I. Lasiecka, {\em Controllability of a viscoelastic Kirchhoff plate} Control and estimation of distributed parameter systems (Vorau, 1988), Internat. Ser. Numer. Math., 91, Birkhäuser, Basel, 1989, 237--247.

%\bibitem{LasT} I. Lasiecka, R. Triggiani, {\em  Exact controllability of the wave equation with Neumann boundary control} Appl. Math. Optim. {\bf 19} (1989),  243–290.


%\bibitem{Ko} R. Kohlrausch, Theorie des elektrischen r\"uckstandes in der leidener flasche, Pogg. Ann. Phys. Chem., 91 (1854), 179.
%
%\bibitem{LPC} G. Lebon, C. Perez-Garcia, J. Casas-Vazquez,
%{\em On the thermodynamic foundations of viscoelasticity}  J. Chem. Phys., {\bf
%88} (1988),  5068--5075.

%\bibitem{L0}  G. Leugering, {\em Exact boundary controllability of an integro-differential equation},
%Appl. Math. Optim., (1987), 223--250.
%
%\bibitem{L}  G. Leugering, {\em Boundary controllability of a viscoelastic string},
%in G. Da Prato and M. Iannelli editors, Volterra integrodifferential equations in Banach spaces and
%applications, Harlow, Essex, Longman Sci. Tech., (1989), 258-270.
%
%\bibitem{LR}
%T. Li, B. Rao, {\em Exact synchronization for a coupled system of wave equations with Dirichlet boundary controls}, Chin. Ann. Math. Ser. B {\bf 34} (2013), 139--160.
%
%\bibitem{Lio1} J.-L. Lions, 
%{\em Exact controllability, stabilization and perturbations for distributed systems}, 
%SIAM Rev. {\bf 30} (1988), 1--68. 
%
%\bibitem{Lio2} J.-L. Lions,
%{\em Contr\^olabilit\'e exacte, perturbations et stabilisation de syst\`emes distribu\'es. Tome 1. Contr\^olabilit\'e exacte}, with appendices
%by E. Zuazua, C. Bardos, G. Lebeau and J. Rauch. Recherches en Math\'ematiques Appliqu\'ees, {\bf 8} ( 1988), Masson, Paris.
%
%\bibitem{Lio3} J.-L. Lions,
%{\em Contr\^olabilit\'e exacte, perturbations et stabilisation de syst\`emes distribu\'es. Tome 2. Perturbations}, Recherches en Math\'ematiques Appliqu\'ees,
%{\bf 9} (1988), Masson, Paris.


%\bibitem{LPS} P. Loreti, L. Pandolfi, D. Sforza,  {\em Boundary controllability and observability of a viscoelastic string}, SIAM J. Control Optim., {\bf 50} (2012),  820--844.


\bibitem{LoretiSforza} P. Loreti, D. Sforza {\em  Exact reachability 
for second-order integro-differential equations. } C. R. Math. Acad. 
Sci. Paris {\bf 347} (2009),  1153--1158.

\bibitem{LoretiSforza1} P. Loreti, D. Sforza, {\em Reachability problems 
for a class of integro-differential equations}.  J. Differential 
Equations  {\bf 248} (2010), 1711--1755.

%\bibitem{LoretiSforza2} P. Loreti D. Sforza, {\em Multidimensional controllability problems with memory}, M. Ruzhansky, J. Wirth (Eds.) ``Modern Aspects of the Theory of Partial Differential Equations"  Vol 216 of ``Operator Theory: Advances and Applications" 261--274, Birkh\"auser/Springer, Basel, 2011.

\bibitem{LoretiSforza3} P. Loreti, D. Sforza, {\em Control problems for weakly coupled systems with memory},  J. Differential 
Equations  {\bf 257} (2014), 1879--1938.

\bibitem{LoretiSforza4} P. Loreti, D. Sforza,
{\em Hidden regularity for wave equations with memory}, Riv. Mat. Univ. Parma, {\bf 7} (2016), 391--405.

\bibitem{LoretiSforza5} P. Loreti, D. Sforza, {\em Viscoelastic aspects of glass relaxation models}, Phys. A {\bf 526} (2019), 120768, 10 pp. 

%\bibitem{LV} P. Loreti, V. Valente, {\em Partial exact controllability for spherical membranes},  SIAM J. Control
%Optim.,  {\bf 35} (1997),  641-653.

\bibitem{MM} J. C. Mauro,  Y. Z. Mauro, 
{\em On the Prony series representation of stretched exponential relaxation},
Phys. A {\bf 506} (2018), 75--87.


\bibitem{MZ} J. C. Mauro, E. D. Zanotto, Two centuries of glass research: Historical trends, current status, and grand Challenges for the Future, Int. J. Appl. Glass Sci. 5 (2014) 313.    

\bibitem{MPVP} J. C. Mauro, C. S. Philip, D. J. Vaughn, M. S. Pambianchi, Glass science in the United States: Current status and future directions, Int. J. Appl. Glass Sci. 5 (2014) 2.




%\bibitem{MN} J. E. Mu\~noz Rivera, M. G. Naso, {\em Exact controllability for hyperbolic thermoelastic systems with large memory}, Adv. Differential Equations, {\bf 9} (2004), 1369--1394. 
%
%\bibitem{NSW} M. Najafi, G. R. Sarhangi, H. Wang {\em  The study of the  stabilizability of the coupled wave 
%equations under various end conditions}, in Proceedings of the 31st IEEE Conference  on Decision and Control, Vol. I (Tucson Arizona 1992), 374--379, IEEE, New York, 1992.

%\bibitem{N} M. Najafi,????

%\bibitem{Pruss} J. Pr\"uss, Evolutionary integral equations and
%applications, Monographs in Mathematics, 87 (1993), Birkh\"auser Verlag,
%Basel.

%\bibitem{RHN}  M. Renardy, W. J. Hrusa, J. A. Nohel,  Mathematical
%problems in viscoelasticity, Pitman Monographs Pure Appl.Math.,  35
%(1988), Longman Sci.
%Tech., Harlow, Essex.

%\bibitem{Re1} M. Renardy, Are viscoelastic flows under control or out of control? Systems Control Lett.,  54 (2005), 1183--1193.

%\bibitem{Ru} D. L. Russell, 
%{\em Controllability and stabilizability theory for linear partial differential equations: recent progress and open questions} 
%SIAM Rev., {\bf 20} (1978), 639--739. 

%\bibitem{SG}  J. J. Skrzypek, A. W. Ganczarski, 
%Constitutive Equations for Isotropic and Anisotropic Linear Viscoelastic Materials,
%Mechanics of Anisotropic Materials, pp. 57--85, J.J. Skrzypek, A.W. Ganczarski (eds), (2015), Springer, New York.

%\bibitem{T1} R. Triggiani, {\em Exact boundary controllability on $L_2(\Omega)\times H^{-1}(\Omega)$ of the wave equation with Dirichlet boundary control acting on a portion of the boundary $\partial\Omega$, and related problems} Appl. Math. Optim. {\bf 18} (1988), 241--277.

%\bibitem{ZZ1}
%X. Zhang, E. Zuazua, {\em Polynomial decay and control of a 1-d model for fluid-structure interaction} C. R. Math. Acad. Sci. Paris, {\bf 336} (2003),
%745--750.
%\bibitem{ZZ2}
%X. Zhang, E. Zuazua, {\em  Polynomial decay and control of a $1-d$ hyperbolic-parabolic coupled system} J. Differential Equations, {\bf 204} (2004),  380--438.

%\bibitem{ZM} Q. Zheng, J. C. Mauro, Variability in the relaxation behavior of glass: Impact of thermal history fluctuations and fragility, J. Chem. Phys. 146 (2017), 074504.
%
%\bibitem{ZYMK} T. Zhou, J. Yan, J. Masuda, T. Kuriyagawa, Investigation on the viscoelasticity of optical glass in ultraprecision lens molding process J. Materials Processing Technology, 209 (2009), 4484--4489.






\end{thebibliography}
\end{document}